\documentclass[12pt,leqno]{amsart}

\usepackage{amssymb,amsmath,amsthm,amssymb,dsfont}
\usepackage{graphicx,tikz}
\usetikzlibrary{arrows.meta}
\usepackage[margin=0.55in]{geometry}
\usepackage{amsmath,amssymb,booktabs,longtable,array,xcolor}
\usepackage{color}
\usepackage{verbatim}
\usepackage[normalem]{ulem}
\usepackage[bookmarksopen,bookmarksdepth=3,colorlinks,citecolor=red,pagebackref,hypertexnames=true]{hyperref}
\usepackage{amsmath,amssymb,booktabs,longtable,array,xcolor}
\usepackage[msc-links,nobysame,non-sorted-cites,initials]{amsrefs}
\usepackage[inline]{enumitem}
\usepackage{mathtools}
\mathtoolsset{showonlyrefs}

\definecolor{darkblue}{RGB}{0,0,160}

\usepackage[lining]{libertine}
\usepackage{cabin}
\usepackage[libertine]{newtxmath}
\usepackage[utf8]{inputenc}
\usepackage[T1]{fontenc}
\usepackage{indentfirst}
\def\Z {{\mathbb Z}}

\numberwithin{equation}{section}

\newtheorem{proposition}[equation]{Proposition}
\newtheorem{theorem}[equation]{Theorem}
\newtheorem{corollary}[equation]{Corollary}
\newtheorem{lemma}[equation]{Lemma}

\theoremstyle{definition}

\newtheorem*{remark*}{Remark}
\newtheorem*{warn*}{A word of warning}

\theoremstyle{plain}

\usepackage{scalerel}
\newcommand{\EE}{\mathbb{E}}

\begin{document}
		
		\raggedbottom
		
		\title{The Bilinear Hilbert Transform Along the Primes}

		\author{Anastasios Fragkos}
		\address{
			Department of Mathematics,
			Georgia Tech \\
			686 Cherry Street, Atlanta, GA 30332-0160 USA}
		\email{anastasiosfragkos@gatech.edu}

		\author{Ben Krause}
		\address{
			Department of Mathematics,
			University of Bristol \\
			Beacon House, Queens Rd, Bristol BS8 1QU}
		\email{ben.krause@bristol.ac.uk}
		
		\author{Michael Lacey}
		\address{
			Department of Mathematics,
			Georgia Tech \\
			686 Cherry Street, Atlanta, GA 30332-0160 USA}
		\email{lacey@math.gatech.edu}
		
		\author{Nazar Miheisi}
		\address{
			Department of Mathematics,
			King's College London \\
			Strand Building, Strand Campus, Strand, London, WC2R 2LS}
		\email{naz.miheisi@gmail.com}

		\author{Yu-Chen Sun}
		\address{
			Department of Mathematics,
			University of Bristol \\
			Beacon House, Queens Rd, Bristol BS8 1QU}
		\email{yuchensun93@163.com}
		
		\date{\today}
		
		\begin{abstract}
			We prove that the discrete Bilinear Hilbert Transform along the primes,
			\begin{align}
				\mathcal{H}_{\mathbb{P}}(f,g)(x) := \sum_{n \neq 0} \frac{\Lambda(|n|)}{n} f(x-n) g(x+n), \qquad \Lambda(n) := \begin{cases} \log p & \text{ if } n = p^k, \ k \geq 1 \\
            0 & \text{ otherwise} \end{cases}
			\end{align}
			extends to a bounded operator from $\ell^p \times \ell^q \to \ell^r$ whenever $p,q \in (1,\infty)$ and $\frac{1}{r}=\frac{1}{p}+\frac{1}{q}$ in the Banach range, $r \geq 1$, and a particular portion of the quasi-Banach range $r \in \left(\frac{2}{3},1\right)$; in the diagonal case, this extends up to the endpoint $(\frac{5}{3},\frac{5}{3},\frac{5}{6})$.  After $U^3$-based reductions, the philosophy of the proof becomes similar to that of Krause-Mirek-Tao (on different bilinear operators in discrete harmonic analysis), namely a ``splitting" of the problem into an arithmetic component, handled by analysis within cyclic subgroups, and an analytic component, handled by analysis along arithmetic progressions. For the arithmetic component, we draw upon previous work of most of the authors on Carleson operators along the primes; the main analytic input are estimates for \emph{UMD space}-valued bilinear singular integrals with modulation invariance, due to Di Plinio-Li-Martikainen-Vuorinen.
		\end{abstract}
		
		\maketitle
		
		\setcounter{footnote}{0} 

		\setcounter{tocdepth}{1}
		
		\section{Introduction}
		
		This paper is a continuation of the authors' program exploring discrete analogues in harmonic analysis along the primes that feature modulation invariance, see \cite{F+,FKT1}, and, most significant to our current context, the Carleson operator along the primes \cite{FKMS26},
		\begin{align}
			\mathcal{C}_{\mathbb{P}}f(x) := \sup_{\lambda \in \mathbb{T}} \big| \sum_{n \neq 0} \frac{\Lambda(|n|)}{n} f(x-n) e^{2 \pi i n \lambda} \big|;
		\end{align}
		our focus here is the Bilinear Hilbert Transform along prime times,
		\begin{align}\label{e:BHTdef}
			\mathcal{H}_{\mathbb{P}}(f,g)(x) := \sum_{n \neq 0} \frac{\Lambda(|n|)}{n} f(x-n) g(x+n),
		\end{align}
		where $\Lambda$ denotes the von Mangoldt function which can be thought as an $L^1$-normalized analogue of the indicator of primes. Our main result is stated below
		\begin{theorem} \label{t:maintheorem}
			Let $\mathcal{R} \subseteq (0,\infty)^{3}$ be defined as the region  \[  \mathcal{R} \coloneqq  \left\{ (p,q,r) \in (0,\infty)^3: \frac{2}{p}+\frac{3}{q}<3,\; \; \frac{3}{p}+\frac{2}{q}<3, \;\; \frac{1}{r}=\frac{1}{p}+\frac{1}{q}, \; p,q>1  \right\} \] then whenever $(p,q,r) \in \mathcal{R}$ we have that
			\[ \| \mathcal{H}_{\mathbb{P}}(f,g) \|_{\ell^r} \leq C_{p,q,r} \|f \|_{\ell^p} \| g \|_{\ell^q}.\]
			In particular, the above estimate holds whenever $r = \frac{5}{6}^+, \ p = q = 2r = \frac{5}{3}^+$.
		\end{theorem}
		\begin{figure}[tbp]
			\centering
			\begin{tikzpicture}[
				x=7.5cm,y=7.5cm,
				font=\small,
				line cap=round,line join=round,
				>=Stealth
				]
				\definecolor{regionblue}{RGB}{44,82,130}
				\definecolor{regionorange}{RGB}{215,125,35}
				
				\coordinate (O) at (0,0);
				\coordinate (A) at (1,0);
				\coordinate (B) at (0,1);
				\coordinate (D) at (0.6,0.6);
				
				\fill[regionblue!12] (O)--(A)--(B)--cycle;
				\fill[regionorange!25] (A)--(D)--(B)--cycle;
				
				\draw[->,gray!65,line width=0.45pt]
				(-0.04,0)--(1.14,0)
				node[right,black] {$a=\frac{1}{p}$};
				\draw[->,gray!65,line width=0.45pt]
				(0,-0.04)--(0,1.14)
				node[above,black] {$b=\frac{1}{q}$};
				
				\draw[black!65,densely dotted,line width=0.8pt]
				(A)--(B);
				\node[rotate=-45,fill=regionblue!12,inner sep=2pt]
				at (0.73,0.22) {$a+b=1\quad(r=1)$};
				
				\draw[regionblue,dashed,line width=0.9pt]
				(O)--(A)--(D)--(B)--cycle;
				
				\foreach \P in {O,A,B,D}
				\draw[regionblue,fill=white,line width=0.8pt]
				(\P) circle[radius=1.8pt];
				
				\node[below left=3pt] at (O) {$(0,0)$};
				\node[below=5pt] at (A) {$(1,0)$};
				\node[left=5pt] at (B) {$(0,1)$};
				\node[above right=4pt] at (D)
				{$\left(\frac{3}{5},\frac{3}{5}\right)$};
				
				\node[rotate=-56.31,text=regionblue]
				at (0.85,0.34) {$3a+2b=3$};
				\node[rotate=-33.69,text=regionblue]
				at (0.30,0.88) {$2a+3b=3$};
				
				\node[text=regionblue,font=\Large]
				at (0.28,0.34) {$\mathcal R$};
				\node at (0.29,0.25) {$r>1$};
				\node[text=regionorange!75!black]
				at (0.53,0.545) {$r<1$};
			\end{tikzpicture}
			\caption{The region $\mathcal R$ in reciprocal coordinates
				$a=\frac{1}{p}$ and $b=\frac{1}{q}$.
				Dashed boundaries are excluded; the dotted $r=1$
				segment is included.}
			\label{fig:region-R}
		\end{figure}
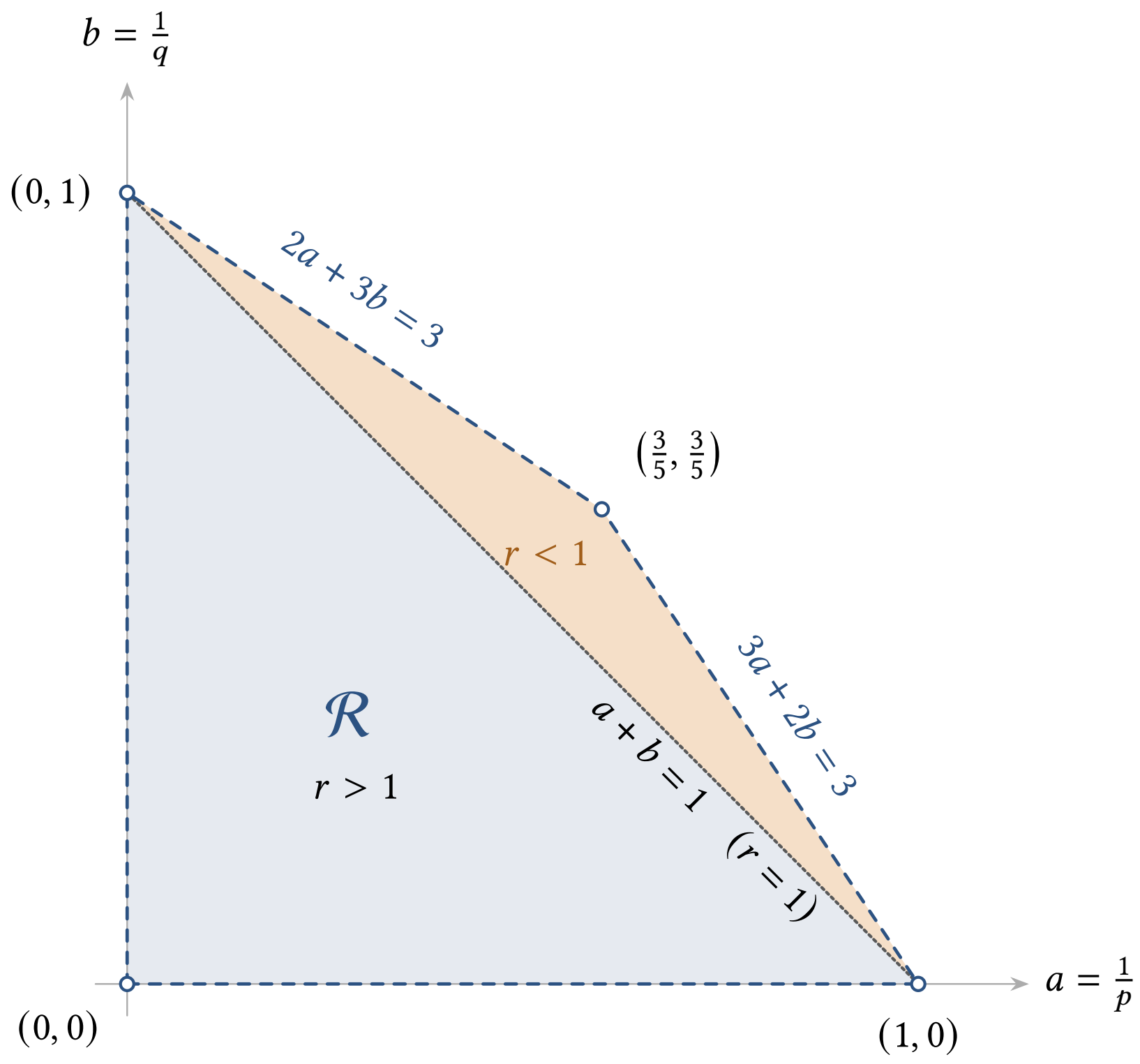
		\subsection{History}
		The Bilinear Hilbert Transform,
		\begin{align}
			H(f,g)(x) :=  \operatorname{p.v.} \int_{\mathbb{R}} f(x-t) g(x+t) \ \frac{dt}{t}
		\end{align}
		was first introduced by Calder\'{o}n in his study of the first commutator of the Cauchy integral along Lipschitz curves, see \cite{MSLipschitz} for a nice discussion. The celebrated articles by the third author and C. Thiele \cite{LT1999,LT2000} established the following norm estimates for $H$.
		\begin{theorem}\label{t:BHT}
			Suppose that $f \in L^p(\mathbb{R}), \ g \in L^q(\mathbb{R})$ and 
			\begin{align}\label{e:fullbht} \frac{1}{p} + \frac{1}{q} = \frac{1}{r} < \frac{3}{2}, \qquad p,q > 1.\end{align}
			Then there exists an absolute constant $0 < C_{p,q,r} < \infty$ so that
			\begin{align}
				\| H(f,g) \|_{L^r(\mathbb{R})} \leq C_{p,q,r} \|f \|_{L^p(\mathbb{R})} \| g \|_{L^q(\mathbb{R})}.
			\end{align}
		\end{theorem}
		This result was later strengthened by the third author \cite{Lacey2000} to incorporate maximal truncations
		\begin{align}
			H(f,g)(x) \longrightarrow H^*(f,g)(x) := \sup_{0<r<R<\infty} \big| \int_{r < |t| < R} f(x-t) g(x+t) \ \frac{dt}{t} \big|
		\end{align}
		under the same integrability conditions; for the remainder of the paper, we will refer to functions in integrability classes satisfying \eqref{e:fullbht} as being in the (full) BHT range.
		
		Standard transference arguments imply sequence-space analogues, which include the bound
		\begin{align}
			\left  \| \sum_{n \neq 0} \frac{f(x-n) g(x+n)}{n} \right\|_{\ell^r} \leq C_{p,q,r} \|f \|_{\ell^p} \| g \|_{\ell^q}
		\end{align}
		in the full BHT range. And, by suitably refining \cite{Lacey2000}, see \cite{BeckerDurcik26,Demeter07,DoOberlinPalsson17}, and applying Calder\'{o}n's transference principle \cite{Calderon1968}, one can establish a pointwise convergence statement, the bilinear analogue of Cotlar's ergodic theorem \cite{Cotlar1955}.
		\begin{theorem}\label{t:BILINEARCOTLAR}
			If $(X,\mu,T)$ is an invertible measure-preserving system on a $\sigma$-finite measure-space, then the ergodic sums
			\begin{align}
				\sum_{1 \leq |n| \leq N} \frac{f(T^n x) g(T^{-n} x)}{n}
			\end{align}
			converge $\mu$-almost everywhere for functions $f,g$ in the full BHT range.
		\end{theorem}
		
		The focus of this paper is the discrete Bilinear Hilbert Transform along prime times, which combines modulation invariance with Radon-type behavior. Indeed, if we conflate the notion of sequence-space density with dimension, the asymptotic
		\begin{align}
			N^{-o(1)} \leq \frac{|\{p \leq N: \text{ prime}\}|}{N} \approx \frac{1}{\log N} = o_{N\to \infty}(1)
		\end{align}
		provides a perspective on $\Lambda(|\cdot|)$ as roughly analogous to an even, singular, Ahlfors-David regular of full Hausdorff dimension. But, we emphasize that no Lebesgue-singular model for $H$ has been studied, so that, in analogy with the results of \cite{FKMS26}, the discrete Radon theory has progressed further than its continuous analogue; we refer the reader to the introduction of \cite{FKMS26} for further discussion. As discussed therein, the study of discrete operators supported on the primes was motivated by questions in pointwise ergodic theory, and in this regard we emphasize that norm estimates on the maximally truncated version of $\mathcal{H}_{\mathbb{P}}$
		\begin{align}
			\left \| \sup_{N} \left| \sum_{1 \leq |n| \leq N} \frac{\Lambda(|n|)}{n} f(x-n) g(x+n) \right| \right \|_{\ell^r} \leq C_{p,q,r} \|f \|_{\ell^p} \| g \|_{\ell^q}, \qquad 0 < C_{p,q,r} < \infty,
		\end{align}
		are essentially sufficient to establish the analogue of Theorem \ref{t:BILINEARCOTLAR} along prime times, after a brief interpolation argument; we hope to address this in future work, and confine our analysis below to Theorem \ref{t:maintheorem}.
		
		\subsection{Proof Strategy}
		At the highest level, the philosophy of proof is a synthesis of \cite{KMT} and \cite{FKMS26}: 
		
		Following $U^3$ reductions, the task at hand is to estimate the bilinear single complexity operators
		\[
		\sum_n K_{Q}(n)w_{Q,B}(n)f(x-n)g(x+n),
		\] associated to a block of square-free integers $B \subset (Q/2,Q]$, 
		with sufficient gain in $Q$, where $K_Q$ is a suitably smooth Calder\'{o}n-Zygmund kernel, and $w_{Q,B}$ is as in \eqref{e:WQB}; so far, this is analogous to \cite{FKMS26}. In that work, ``variable coefficient'' multi-frequency analysis was employed to reduce matters to obtaining gain at the single-scale level in the cyclic-subgroup setting, at which point Hilbert-space orthogonality and elementary number-theoretic methods were sufficient to conclude.
		
		In our current context, while the periodic analysis of \cite{FKMS26} adapts readily to the bilinear setting, we are unable to directly ``lift'' the relevant periodic estimates. To get around this point, we adopt the perspective of \cite{KMT} and view our operator from a ``tensor'' perspective, as an analytic operator applied to Schwartz functions taking values in cyclic subgroups. With this perspective in mind, work of \cite{DPLMV} on modulation invariant operators taking values in so-called \emph{intermediate UMD spaces} -- of which 
		\[ \ell^p(\mathbb{Z}/N), \qquad  1 < p < \infty, \ N \geq 1, \]
		form examples -- interposes, and we are able to propagate the periodic gain, after taking into account \emph{aliasing} issues involved in passing from the Euclidean- to sequence-space setting.
		
		\medskip
		
		The sections are organized as follows, with the final interpolation argument deferred until all estimates are available:
		\begin{itemize}
			\item Section \ref{s:Prelim} records preliminary material: basic notational conventions and Gowers information;
			\item Section \ref{s:red} applies $U^3$-based arguments to reduce matters to our key Proposition \ref{p:key0};
			\item Section \ref{s:per} develops our analysis in the periodic setting; and
			\item Section \ref{s:transference} transfers the periodic gain through the vector-valued theorem of \cite{DPLMV}; and
			\item Section~\ref{s:final} completes the proof by interpolation.
		\end{itemize}

		\section{Preliminaries}\label{s:Prelim}
		
		\subsection{General Notation}
		Throughout, we let
		\begin{align}
			e(t) := e^{2 \pi i t}, \; \; \; e_N(t) := e(t/N)
		\end{align}
		denote the complex exponential, while for sequences on $\mathbb Z$ and
		functions on $\mathbb R$ our Fourier conventions are
		\[
		\widehat f(\beta)=\sum_{n\in\mathbb Z}f(n)e(-n\beta),\quad\beta\in\mathbb T,
		\qquad
		\widehat F(\xi)=\int_{\mathbb R}F(x)e(-x\xi)\,dx, \quad \xi \in \mathbb{R};
		\]
		all other Fourier conventions will be specified locally.
		Both inverse transforms use the opposite sign in the exponential;
		the torus has normalized Lebesgue measure. We 
		identify 
		\[ \mathbb{Z}/N := \{ 0,1,\dots, N-1\} \mod N,\]
		use the standard notation for $L^1$-normalized dilations,
		\begin{align}
			\varphi_N(n) := \frac{1}{N} \varphi(\frac{n}{N}),
		\end{align}
		and let
		\begin{align}
			M_{\text{HL}} f(x):=  \sup_{I \ni x} \, \frac{1}{|I|} \sum_{n \in I} |f(n)|  
		\end{align}
		denote the discrete Hardy-Littlewood maximal function; here and throughout, we use $I$ to denote discrete intervals. We let
		\begin{align}\label{e:bimax}
			&B^*(F,G)(x) := \sup_L B_L(F,G), \qquad B_s^*(F,G)(x) := \big(B^*(|F|^s,|G|^s)(x)\big)^{1/s} \\ & B_L(F,G)(x):=\frac{1}{2L+1}\sum_{|n|\leq L}|F(x-n)G(x+n)|
		\end{align}
		denote Lacey's bilinear maximal function and its power variant, where the former is bounded in the full BHT range.
		
		We set $w(n) := \Lambda(|n|)$ for $n\ne0$ and $w(0):=0$, where $\Lambda$ is the von Mangoldt function:
		$\Lambda(n)=\log p$ if $n$ is a positive power of the prime $p$, and $\Lambda(n)=0$ otherwise.
		
		We collect
		\begin{align}\label{e:GammaQ}
			\Gamma_Q := \{ a/q \text{ reduced} : Q/2 < q \leq Q \},
		\end{align}
		and define
		\begin{align}\label{e:WQ}
			\Lambda_Q(n) := \sum_{a/q \in \Gamma_Q} \frac{\mu(q)}{\phi(q)} e(na/q ) = \sum_{Q/2 < q \leq Q} \frac{\mu(q)}{\phi(q)} c_q(n), \qquad w_Q(n) := \Lambda_Q(|n|)
		\end{align}
		where
		\begin{align}\label{e:ram}
			c_q(n) = \sum_{a/q : \ (a,q) =1} e(na/q )
		\end{align}
		denotes a Ramanujan sum, and $\mu,\phi$ denote the M\"{o}bius and totient functions, respectively; we will let $d_{\text{div}}$ denote the divisor function, well-known to grow sub-polynomially; we will often make use of the lower bound
		\begin{align}\label{e:totientbound}
			\phi(q) \gtrsim \frac{q}{\log \log (10 + q)}.
		\end{align}
		
		For (square-free) subsets $B \subset (Q/2,Q]$
		we define
		\begin{align}\label{e:lcmQ}
			\mathcal{Q}_B := \text{lcm}( q \in B )
		\end{align}
		and analogously define
		\begin{align}\label{e:WQB}
			\Lambda_{Q,B}(n) := \sum_{q \in B} \frac{\mu(q)}{\phi(q)} c_q(n), \qquad w_{Q,B}(n) := \Lambda_{Q,B}(|n|).
		\end{align}
		
		Choose an even function 
		\begin{align}\label{e:chi} 
			\chi\in\mathcal C_c^\infty((-1,1))
		\end{align}
		which equals one on $[-1/2,1/2]$, and set
		\[
		\psi(t):=\frac{\chi(t)-\chi(2t)}{t},\qquad \psi(0):=0.
		\]
		This is a smooth odd function supported where $1/4\leq|t|<1$, so that
		\begin{align}
			\sum_{k\geq1}\psi_k(t)=\frac1t,\quad |t|\geq1,
			\qquad \psi_k(t):=2^{-k}\psi(t2^{-k}).
		\end{align}
		
		Whenever the exponents $(p,q,r)$ occur below, we assume $1<p,q<\infty$ and 
		\[ \frac{1}{p} + \frac{1}{q} = \frac{1}{r}.\]
		
		\subsection{Asymptotic Notation}\label{sss:O}
		We use the modified Vinogradov notation, where we write $X \lesssim Y$ or $Y \gtrsim X$ to denote
		the estimate $X \leq CY$ for an absolute constant $C$ and $X, Y \geq 0.$  If we need $C$ to depend on a
		parameter, we shall indicate this by subscripts, so for instance $X \lesssim_p Y$ denotes the estimate $X \leq C_p Y$ for some $C_p$ depending on $p$. We use $X \approx Y$ as shorthand for $Y \lesssim X \lesssim Y$, while the notation $X \ll Y$ or $Y \gg X$ indicates that the implicit constant in the $\lesssim$ notation is extremely large, and analogously $X \ll_p Y$ and $Y \gg_p X$.
		
		We also make use of big-Oh and little-oh notation: we let $O(Y)$  denote a quantity that is $\lesssim Y$, and similarly
		$O_p(Y )$ will denote a quantity that is $\lesssim_p Y$; we let $o_{t \to a}(Y)$
		denote a quantity whose quotient with $Y$ tends to zero as $t \to a$ (possibly $\infty$), and
		$o_{t \to a;p}(Y)$
		denote a quantity whose quotient with $Y$ tends to zero as $t \to a$ at a rate depending on $p$. When clear from context, we will suppress the $t \to a$ subscript.

		\subsection{Gowers Norms}
		\label{sub:gowers_norms}
		For a finitely supported function $f\colon\mathbb Z\to\mathbb C$, define its conjugation differences by 
		\begin{align}\label{e:triangle}
			\Delta_h f(x) \coloneqq  f(x) \overline{f(x+h)}, 
			\qquad x,h\in \mathbb{Z}, 
		\end{align}
		so that 
		\begin{align} \label{e;double}
			\Bigl\lvert \sum_x f(x)\Bigr\rvert ^2 
			= 
			\sum_{x,h} \Delta_h f(x). 
		\end{align}
		The higher order conjugation-difference operator is inductively defined to be 
		\begin{align}
			\Delta_{h_1,\dots,h_s}f(x) \coloneqq  \Delta_{h_s}(\Delta_{h_1,\dots,h_{s-1}} f)(x), 
			\qquad x,h_1, \ldots, h_s\in \mathbb{Z} . 
		\end{align}
		Then, the $s$th order Gowers norm is 
		\begin{align} \label{e;GowersDef}
			\| f \|_{U^s(\mathbb{Z})}^{2^s} \coloneqq  
			\sum_{x,h_1,\dots,h_s} \Delta_{h_1,\dots,h_s}f(x). 
		\end{align}
		For $s=1$, this is a semi-norm, while higher orders are norms. In particular, for $s=2$, we have 
		\begin{align}  \label{e;U2-ell4}
			\lVert f \rVert_{U^2(\mathbb{Z} )} ^{4} 
			= \int_{\mathbb{T} } \lvert  \widehat{f} (\beta )\rvert ^{4}\;d \beta;
		\end{align}  
		we recall the inductive relationship between norms given by 
		\begin{align}\label{e:foliate}
			\| f \|_{U^{s+1}(\mathbb{Z})}^{2^{s+1}} = \sum_{h_1,\dots,h_{s-1}} \| \Delta_{h_1,\dots,h_{s-1}} f \|_{U^2(\mathbb{Z})}^4.
		\end{align}

		When $N$ is a positive integer, we write $[N]=\{1,2,\dots,N\}$ and denote the usual expectation by 
		\begin{align}
			\mathbb{E}_{n \in [N]} f(n) := \frac{1}{N} \sum_{n \in [N]} f(n)\end{align}
		and more generally, for finite $X \subset \mathbb{Z}$, define
		\begin{align}
			\mathbb{E}_{n \in X} f(n) := \frac{1}{|X|} \sum_{n \in X} f(n) \end{align}
		and the norms
		\begin{align}
			\| f \|_{L^p(X)}^p := \frac{1}{|X|} \sum_{n \in X} |f(n)|^p,
		\end{align}
		and inner products
		\begin{align}\label{e:Innerprod}
			\langle f,g \rangle_X := \frac{1}{|X|} \sum_{n \in X} f(n) \overline{g(n)}.
		\end{align}
		Additionally,  we introduce the periodic Fourier transform
		\begin{align}
			\mathcal{F}_N f(\xi) := \mathbb{E}_{r \in [N]} f(r) e_N(-\xi r) : \mathbb{Z}/N \to \mathbb{C},
		\end{align}
		where we recall $e_N(t) := e(t/N)$.
		
		We define normalized $U^s$ norms as follows: for $f \colon [N] \to \mathbb C$, 
		set 
		\begin{equation}
			\lVert f \rVert_{U^s([N]) } \coloneqq  
			\frac{\| f \|_{U^s(\mathbb{Z})}}{\lVert \mathbf{1}_{[N]} \rVert_{U^s(\mathbb Z) }} \approx \frac{\| f \|_{U^s(\mathbb{Z})}}{N^{\frac{s+1}{2^s}}}.
		\end{equation}
		
		We make repeated use of the $s=3$ case of the following estimate, due to Eisner-Tao \cite{ET}, relating $U^s$ and $\ell^p$ norms for appropriate $p$ (determined by scaling).
		
		\begin{lemma}\label{l:ET}
			If $s\geq1$ and $p_s:=\frac{2^s}{s+1}$, then
			\begin{align}
				\| f \|_{U^s([N])} \lesssim_s \| f \|_{L^{p_s}([N])}.
			\end{align}
		\end{lemma}
		
		We will also make use of the following $U^3$ estimate; for completeness we include a proof.
		\begin{lemma}\label{l:U3maxmod}
			For $W:[N]\to\mathbb C$ and finitely supported $f,g:\mathbb Z\to\mathbb C$,
			\begin{align}\label{e:U3global}
				\left\|\mathbb E_{n\in[N]}W(n)f(x-n)g(x+n)\right\|_{\ell^2_x}
				\lesssim \|W\|_{U^3([N])}\|f\|_{\ell^4}\|g\|_{\ell^4}.
			\end{align}
			The same bound holds with the average replaced by
			$\sum_n\phi_N(n)W(n)f(x-n)g(x+n)$, whenever
			$\phi\in\mathcal C_c^2([0,1])$ and
			$\|\phi\|_{\mathcal C^2(\mathbb R)}\leq1$.
		\end{lemma}
		\begin{proof}
			Let $L$ denote the left-hand side of \eqref{e:U3global}, and extend
			$W$ by zero outside $[N]$.
			Expanding the square and translating the spatial variable, we obtain
			\[
			L^2=N^{-2}\sum_{h,n,x}
			\Delta_hW(n)\Delta_{-h}f(x)\Delta_hg(x+2n),
			\]
			where only $|h|<N$ contributes.
			Since $\sum_{|h|<N,x}|\Delta_{-h}f(x)|^2\lesssim
			N\|f\|_{\ell^4}^4$, Cauchy--Schwarz followed by expansion of the
			remaining square implies
			\begin{align*}
				L^4
				&\lesssim N^{-3}\|f\|_{\ell^4}^4
				\sum_{|h|<N,x}\left|\sum_n\Delta_hW(n)\Delta_hg(x+2n)\right|^2\\
				&\lesssim N^{-3}\|f\|_{\ell^4}^4
				\sum_{|h|,|k|<N}
				\left(\sum_n\Delta_{h,k}W(n)\right)
				\left(\sum_x\Delta_{h,2k}g(x)\right).
			\end{align*}
			Applying Cauchy--Schwarz in $(h,k)$, and using H\"older to bound
			$|\sum_x\Delta_{h,2k}g(x)|\leq\|g\|_{\ell^4}^4$, we conclude that
			\begin{align*}
				L^8
				&\lesssim N^{-6}\|f\|_{\ell^4}^8
				\sum_{h,k}\left|\sum_n\Delta_{h,k}W(n)\right|^2
				\sum_{|h|,|k|<N}\left|\sum_x\Delta_{h,2k}g(x)\right|^2\\
				&\lesssim N^{-4}\|W\|_{U^3(\mathbb Z)}^8
				\|f\|_{\ell^4}^8\|g\|_{\ell^4}^8.
			\end{align*}
			The normalization of $U^3([N])$ proves the first assertion, while
			Fourier inversion of $\phi$ proves the second because
			$\|\widehat\phi\|_{L^1(\mathbb R)}\lesssim\|\phi\|_{\mathcal C^2(\mathbb{R})}$
			and linear modulations preserve the $U^3$ norm.
		\end{proof}

		\subsection{Single Scale Estimates}
		Before turning to the main argument, we record the estimates that will
		control the first group of scales and the intermediate groups.
		\begin{lemma}\label{l:singlescalecomplexitybetter}
			If $k\geq b_Q$, with $b_Q$ defined in \eqref{e:bQ} below, then
			\begin{align}\label{e:single-L1}
				\left\|\sum_n w_Q(n)\psi_k(n)f(x-n)g(x+n)\right\|_{\ell^1}
				&\lesssim Q^{-1/4+o(1)}\|f\|_{\ell^2}\|g\|_{\ell^2}.
			\end{align}
			Whenever $p,q>1$ and $1/r=1/p+1/q$, we also have
			\begin{align}\label{e:single-rough}
				\left\|\sum_n w_Q(n)\psi_k(n)f(x-n)g(x+n)\right\|_{\ell^r}
				&\lesssim_{p,q}Q^{o(1)}\|f\|_{\ell^p}\|g\|_{\ell^q},
			\end{align}
			so the second estimate includes the range $1/2<r\leq 1$.
		\end{lemma}
		\begin{proof}
			Whenever $u,v\geq1$ and $1/t=1/u+1/v$, we have
			\begin{equation}\label{e:fixed-average}
				\|B_L(F,G)\|_{\ell^t}\lesssim_{u,v}
				\|F\|_{\ell^u}\|G\|_{\ell^v}.
			\end{equation}
			For $t\geq1$ this follows from Minkowski and H\"older, while for
			$1/2\leq t<1$ we partition the output into intervals $I$ of length
			comparable to $L$ and restrict both inputs to $3I$.
			At $(1,1,2)$ we have
			\[
			\|B_L(F\mathbf1_{3I},G\mathbf1_{3I})\|_{\ell^{1/2}(I)}
			\leq  \| F \|_{\ell^{1}(3I)} \| G \|_{\ell^1(3I)}
			\]
			and since every tuple of interest is at the convex hull of $(1,1,2), (1,0,1), (0,1,1)$ proves \eqref{e:fixed-average},
			since the enlarged intervals have bounded overlap. With this in mind, we appeal to convexity: fix $1<s<\min\{p,q\}$, and apply the moment estimates of \cite[Lemma 4.6]{KMTT25} to $\Lambda_{\leq Q}-\Lambda_{\leq Q/2}$; H\"older in $n$ implies
			\[
			\left|\sum_n w_Q(n)\psi_k(n)f(x-n)g(x+n)\right|
			\lesssim_s Q^{o(1)}
			\big(B_{2^k}(|f|^s,|g|^s)(x)\big)^{1/s},
			\]
			which proves \eqref{e:single-rough} by \eqref{e:fixed-average}.
			By \cite[Proposition 2.1]{F+} and Lemma~\ref{l:U3maxmod}, the
			$(4,4,2)$ norm is $O(Q^{-3/8+o(1)})$, provided $2^k \geq Q^{20}$.
			Interpolating with \eqref{e:single-rough} at $(u,u,u/2)$, where
			$u>1$ is arbitrarily close to one, proves \eqref{e:single-L1};
			the finitely many smaller values of $Q$ are absorbed into the constant.
		\end{proof}
		With these preliminaries in hand, we turn to the main body of our work.
		\section{Reductions and Approach}\label{s:red}
		Since all sums in $Q$ will be dyadic, we take $Q\in\{1,2,4,\ldots\}$
		and retain the convention $w_1=1$.
		For a finite upper truncation $J$, define
		\[
		H^J(f,g)(x):=\sum_{k=1}^J\sum_n
		\psi_k(n)w(n)f(x-n)g(x+n),
		\]
		so that it suffices to establish bounds independent of $J$.
		Following \cite{F+,FKMS26}, we use the Heath--Brown approximants
		\[
		\Lambda_{\leq M}:=\sum_{Q\leq M}\Lambda_Q,\qquad
		w_{\leq M}:=\Lambda_{\leq M}(|\cdot|),\qquad
		\mathcal E_k:=w-w_{\leq M_k},
		\]
		where
		\begin{align}\label{e:bQ}
			M_k:=\exp((\log 2^k)^{1/10}),\qquad
			b_Q:=\min\{k\geq1:Q\leq M_k\}
			=\max\left\{1,\left\lceil\frac{(\log Q)^{10}}{\log2}\right\rceil\right\}.
		\end{align}

        For the positive and negative halves of the scale $2^k$, the estimates
		of \cite[Propositions 3.1]{F+} and \cite[Proposition 2.1]{F+} respectively read
		\begin{align}\label{e:arithmeticU3}
			\|\mathcal E_k\|_{U^3([2^k])}&\lesssim_A k^{-A},\notag\\
			\|\Lambda_Q\|_{U^3([N])}&\lesssim_\delta Q^{-3/8+\delta},
			\qquad N \geq Q^{20},\ Q\geq2,
		\end{align}
		where the first estimate is used on each half after reflection when needed. We emphasize that the proofs of these highly non-trivial estimates require substantial arithmetic input and essentially encode a large portion of the arithmetic component of our problem. For example, the proof of the \cite[Proposition 3.1]{F+} invokes a quasi-polynomial inverse theorem for the $U^3$ norm.
		Applying Lemma~\ref{l:U3maxmod} on these halves, we obtain
		\begin{align}\label{e:errL2}
			\left\|\sum_n\psi_k(n)\mathcal E_k(n)f(x-n)g(x+n)\right\|_{\ell^2}
			\lesssim_A k^{-A}\|f\|_{\ell^4}\|g\|_{\ell^4}.
		\end{align}
		Since $\Lambda(n)\lesssim k$ on the support of $\psi_k$, the moment
		bounds of \cite[Lemma 4.6]{KMTT25}, followed by the proof of
		\eqref{e:single-rough}, also imply
		\begin{align}\label{e:errrough}
			\left\|\sum_n\psi_k(n)\mathcal E_k(n)f(x-n)g(x+n)\right\|_{\ell^r}
			\lesssim_{p,q} k^{C_{p,q}}\|f\|_{\ell^p}\|g\|_{\ell^q},
		\end{align}
		whenever $p,q>1$ and $1/r=1/p+1/q$.
		Because every such reciprocal tuple is an interior point of the input
		range, interpolation with \eqref{e:errL2}, with $A$ sufficiently large,
		yields
		\begin{align}\label{e:error-summable}
			\left\|\sum_n\psi_k(n)\mathcal E_k(n)f(x-n)g(x+n)\right\|_{\ell^r}
			\lesssim_{p,q} k^{-100}\|f\|_{\ell^p}\|g\|_{\ell^q}.
		\end{align}
		
		Next, for $0<\tau<1$, let
		\begin{align*}
			H^J_{Q,\mathrm{moderate};\tau}(f,g)
			&:=\sum_{\substack{b_Q\leq k<Q^\tau\\k\leq J}}
			\sum_n\psi_k(n)w_Q(n)f(x-n)g(x+n),\\
			U^J_{Q;\tau}(f,g)
			&:=\sum_{\substack{Q^\tau\leq k\leq J}}
			\sum_n\psi_k(n)w_Q(n)f(x-n)g(x+n),
		\end{align*}
		where we use the second definition only for $Q$ large enough that
		$b_Q\leq Q^\tau$.
		If $\mathbf C$ is chosen so that this inequality holds for every
		$Q>\mathbf C$, the exact decomposition is
		\begin{equation}\label{e:exactred}
			\begin{split}
				H^J(f,g)
				&=\sum_{k=1}^J\sum_n\psi_k(n)\mathcal E_k(n)f(x-n)g(x+n)\\
				&\quad+\sum_{Q\leq\mathbf C}\sum_{b_Q\leq k\leq J}
				\sum_n\psi_k(n)w_Q(n)f(x-n)g(x+n)\\
				&\quad+\sum_{Q>\mathbf C}
				\big(H^J_{Q,\mathrm{moderate};\tau}(f,g)+U^J_{Q;\tau}(f,g)\big).
			\end{split}
		\end{equation}
		We suppress the $J$ in the above operators; all subsequent
		estimates remain uniform in the upper truncation parameter $J$. 
		
		We record the progression estimate with both output ranges, since it
		will be used repeatedly in the final interpolation argument.
		If $\mathcal I$ is any finite interval of positive integers, put
		$K_{\mathcal I}:=\sum_{k\in\mathcal I}\psi_k$, so that M\"obius inversion
		expresses its weighted operator as
		\[
		\sum_n K_{\mathcal I}(n)w_Q(n)f(x-n)g(x+n)
		=\sum_{Q/2<q\leq Q}\frac{\mu(q)}{\phi(q)}
		\sum_{d\mid q}\mu(q/d)
		\sum_m dK_{\mathcal I}(dm)f(x-dm)g(x+dm).
		\]
		The innermost operators are uniformly bounded in the full BHT range,
		because the identity
		$tK_{\mathcal I}(t)=\chi(2^{-\max\mathcal I}t)
		-\chi(2^{1-\min\mathcal I}t)$, see \eqref{e:chi},
		represents them as smooth truncations along progressions.
		Using the triangle inequality when $r\geq1$, and the quasi-triangle
		inequality when $r<1$, we conclude that
		\begin{equation}\label{e:progression}
			\left\|\sum_n K_{\mathcal I}(n)w_Q(n)f(x-n)g(x+n)\right\|_{\ell^r}
			\lesssim_{p,q}Q^{\max\{1/r-1,0\}+o(1)}
			\|f\|_{\ell^p}\|g\|_{\ell^q}
		\end{equation}
		in the BHT range: the coefficient sums are bounded by
		\[ \sum_{Q/2<q\leq Q} \frac{d_{\mathrm{div}}(q)}{\phi(q)^s}
		\lesssim Q^{1-s+o(1)}, \qquad s:=\min\{r,1\}.\]
		This controls the bounded-complexity part of \eqref{e:exactred}, while
		the remaining terms require decay in $Q$.
		
		Fix $0<\epsilon<\tau/1000$, and choose a real even Schwartz function
		$\varphi$ with $\widehat\varphi(0)=1$ whose Fourier transform has
		sufficiently small fixed support.
		With $a_Q:=\lceil Q^\tau\rceil$ and $N(Q):=\exp(Q^{\tau-\epsilon/10})$,
		define
		\begin{align}\label{e:KQ}
			K_{Q;\tau}&:=\varphi_{N(Q)}*
			\Big(\sum_{a_Q\leq k\leq J}\psi_k\Big),\notag\\
			\mathcal H_{Q;\tau}(f,g)(x)
			&:=\sum_n K_{Q;\tau}(n)w_Q(n)f(x-n)g(x+n),
		\end{align}
		where convolution is taken on $\mathbb R$ before restriction to
		$\mathbb Z$.
		Because 
		\[ \widehat{\varphi}(0) = \int_{\mathbb{R}} \varphi(t) \ dt =1,\]
		the fundamental theorem of calculus and
		Schwartz decay imply
		\[
		|\varphi_{N(Q)}*\psi_k(n)-\psi_k(n)|
		\lesssim\frac{N(Q)}{2^k}\,2^{-k}(1+|n|/2^k)^{-10},
		\qquad k\geq a_Q,
		\]
		and since $\|\Lambda_Q\|_\infty\leq Q$ and
		$QN(Q)2^{-a_Q}\lesssim Q^{-100}$ for $Q$ sufficiently large, summation over $k$ implies
		\begin{align}\label{e:smootherror}
			|U_{Q;\tau}(f,g)-\mathcal H_{Q;\tau}(f,g)|
			\lesssim Q^{-100}B^*(f,g).
		\end{align}
		The estimate below is the remaining input, whose proof occupies
		Sections~\ref{s:per} and~\ref{s:transference}.
		\begin{proposition}\label{p:key0}
			For every $0<\tau<1$, one may choose the smoothing
			parameter $\epsilon>0$ sufficiently small so that
			\begin{align}\label{e:longgoal}
				\|\mathcal H_{Q;\tau}(f,g)\|_{\ell^1}
				\lesssim_{\tau}Q^{-\tau/2+o(1)}
				\|f\|_{\ell^2}\|g\|_{\ell^2}.
			\end{align}
			Consequently, the unsmoothed operators satisfy
			\begin{align}\label{e:long-unsmoothed}
				\|U_{Q;\tau}(f,g)\|_{\ell^1}
				\lesssim_{\tau}Q^{-\tau/2+o(1)}
				\|f\|_{\ell^2}\|g\|_{\ell^2}.
			\end{align}
		\end{proposition}
		The final section deduces Theorem~\ref{t:maintheorem} from this proposition; we turn to the proof. 
		
		\section{Periodic Analysis}\label{s:per}
		Following the approach of \cite{FKMS26}, we partition the \emph{square-free} integers inside $(Q/2,Q]$ into disjoint subsets $\{ B \}$ with
		\begin{align}\label{e:kappa}
			|B| \approx Q^\kappa, \; \; \; \kappa =\tau  - \epsilon;
		\end{align}
		A last smaller block may be merged with the preceding one; the finitely
		many small $Q$ are already covered by the rough estimates. Possibly after decreasing $\epsilon$, we will assume that $\label{e:epskappa} 0 < \epsilon \leq \frac{\kappa}{100}.$
		
		We record the elementary estimate  $\mathcal{Q}_B:= \operatorname{lcm}\left(q \in B \right) \leq \exp(100 \log Q \cdot Q^\kappa).$
		
		The aforementioned decomposition induces a further splitting of $\mathcal{H}_Q$ into the superposition of  $\left\{\mathcal{H}_{Q,B} : B\right\}$ defined by the equation below, with $w_{Q,B}$ defined in \eqref{e:WQB}. 
		\begin{align}
			\mathcal{H}_{Q;\tau,B}(f,g) :=  \sum_n K_{Q;\tau}(n) w_{Q,B}(n) f(x-n) g(x+n).
		\end{align}
		
		For the remainder of this section, we will fix a single $B$, and regarding $Q$ as fixed, we abbreviate
		\begin{align}
			\Lambda_{Q,B} \longrightarrow \Lambda_B, \qquad w_{Q,B} \longrightarrow w_B;
		\end{align}
		for notational ease, set $   N := \mathcal{Q}_B,$ and regard this parameter as fixed. We abbreviate $Y_N := \mathbb{Z}/N,$ and define
		\begin{align}
			A_B(F,G)(y) := \mathbb{E}_{r \in Y_N} \Lambda_B(r) F(y-r) G(y+r)
		\end{align}
		and the corresponding dual form,
		\begin{align}
			\Pi_B(F,G,H) := \mathbb{E}_{y \in Y_N} A_B(F,G)(y) \overline{H(y)}.    
		\end{align}
		
		For $ p \geq 1$, define
		\begin{align}
			d_B(p) := \| A_B \|_{L^2(Y_N) \times L^2(Y_N) \to L^p(Y_N)}.
		\end{align}
		
		We establish the requisite arithmetic gain in the following proposition; subsequent analytic work will focus on appropriately propagating this gain.
		\begin{proposition}\label{p:interpsetup}
			For each $1 \leq p \leq 2$, we may bound
			\[ d_B(p) \lesssim |B|^{1/2} Q^{-\frac{1}{p}+o(1)}, \]
			where the $o(1)$ is uniform in $B$.
		\end{proposition}
		The $p=2$ case is particularly straightforward.
		\begin{proof}[Proof of the $p=2$ Case]
			Since Cauchy--Schwarz separates the two factors below, we have
			\begin{align*}
				|\Pi_B(F,G,H)|
				&\leq\big(\mathbb E_{y,t}|\Lambda_B(t)|^2|F(y-t)|^2\big)^{1/2}
				\big(\mathbb E_{y,t}|G(y+t)|^2|H(y)|^2\big)^{1/2}\\
				&=\|\Lambda_B\|_{L^2(Y_N)}
				\|F\|_{L^2(Y_N)}\|G\|_{L^2(Y_N)}\|H\|_{L^2(Y_N)}.
			\end{align*}
			The characters occurring in distinct reduced fractions are orthogonal
			on $Y_N$, because each denominator divides $N$, and hence
			\[
			\|\Lambda_B\|_{L^2(Y_N)}^2
			=\sum_{q\in B}\frac1{\phi(q)}\lesssim |B|Q^{-1+o(1)},
			\]
			which proves the claimed operator bound by duality.
		\end{proof}
		By interpolation, Proposition \ref{p:interpsetup} will follow from the bound 
		\[ d_B(1) \leq |B|^{1/2} Q^{-1 + o(1)}. \]
		To prove this result, we argue in analogy with \cite{FKMS26}, passing from our weighted periodic operator, $A_B$, to simpler, unweighted ones; we appeal to convexity to do so. 
		
		Specifically, with $d | q$ for some $q \in B$, set $D := \frac{N}{d}$, and define
		\begin{align}
			A_d (F,G)(y):= \mathbb{E}_{l \in Y_D} F(y-dl) G(y+dl).
		\end{align}
	We can express
		
		\begin{align}\label{eq:divisorsumdecomp}
			A_B(F,G) = \sum_{d |N} \mu(d) \beta_d A_d(F,G), \quad \beta_d := \sum_{q \in B, \ d | q} \frac{1}{\phi(q)},
		\end{align}
		indeed, the identity $c_q(r)=\sum_{d\mid(q,r)}d\mu(q/d)$ and square-free nature of $q \in B$ give
		\[
		\Lambda_B(r)=\sum_{d\mid N}d\mu(d)\beta_d\mathbf{1}_{d\mid r},
		\qquad
		\mathbb{E}_{r\in Y_N}\mathbf{1}_{d\mid r}U(r)
		=d^{-1}\mathbb{E}_{l\in Y_D}U(dl),
		\]
		so the factors $d$ cancel. We have
		\begin{align}
			a_B := \sum_{d | N} \beta_d = \sum_{q\in B}\frac{d_{\text{div}}(q)}{\phi(q)}\lesssim |B| Q^{o(1)-1}.
		\end{align}
		
		As in \cite{FKMS26}, we focus on $\{A_d\}$; fixing $d \mid N$ and noting that $(d,D)=1$, we have the unique decomposition
		\[ x = a + ds, \; \; \; 0\leq a<d, \ s \in \mathbb{Z}/D.\]
		We define $\label{e:fF} f_a(s) := f(a + ds).$ Compute,
		\begin{align}
			\mathcal{F}_D(f_a)(\xi) := \mathbb{E}_{s\in[D]} f_a(s) e_D(-\xi s), \; \; \; \xi \in \mathbb{Z}/D;
		\end{align}
		introduce the fiber-wise maximal function
		\begin{align}
			M_df(a+ds) = \sup_{\xi \in \mathbb{Z}/D} |\mathcal{F}_D(f_a)(\xi)|.
		\end{align}
		\begin{lemma}\label{l:unweighted-fibre}
			Suppose $F,G : Y_N \to \mathbb{C}$; then
			\begin{align}
				\| A_d(F,G) \|_{L^1(Y_N)} \lesssim \| M_d F \|_{L^2(Y_N)} \|G \|_{L^2(Y_N)}.
			\end{align}
		\end{lemma}
		\begin{proof}
			By Fourier inversion, we have
			\begin{align}
				A_d(F,G)(a+ds) = \sum_{\xi \in Y_D} \mathcal{F}_D(F_a)(\xi) \mathcal{F}_D(G_a)(\xi) e( 2\xi s/D);
			\end{align}
			since $\xi \mapsto 2 \xi$ has at most two inverse images, 
			we may bound \[ \EE_{s \in Y_D}|A_d(F,G)(a+ds)|^2 \lesssim \sum_{\xi \in Y_D} |\mathcal{F}_D(F_a)(\xi) \mathcal{F}_D(G_a)(\xi) |^2 \leq (M_dF(a))^2 \| G_a \|_{L^2(Y_D)}^2. \] Hence, since $\| A_d(F,G)\|_{L^1(Y_N)} \leq \| A_d(F,G)(a+ds) \|_{L_{a}^1(Y_d)L_{s}^2(Y_D)} $ we have that \[ \| A_d(F,G)\|_{L^1(Y_N)} \lesssim \EE_{a \in Y_d} M_dF(a) \| G_a \|_{L^2(Y_D)} \leq \| M_dF \|_{L^2(Y_d)} \| G \|_{L^2(Y_N)}=\| M_dF \|_{L^2(Y_N)} \| G \|_{L^2(Y_N)} \]  and the proof is complete.
		\end{proof}
		We continue our analysis of $\{ A_d \}$:
		
		\medskip
		
		Assume \(\beta_d>0\), so \(d\mid q_0\) for some \(q_0\in B\). We emphasize that since \(N\) is square-free and \(N=dD\), $d$ and $D$ are coprime. For \(u\in Y_d\) and \(\rho\in Y_D\), define
		\begin{equation}\label{eq:Kdrho-new}
			\widetilde K_{d,\rho}(u)
			:=\mathbb{E}_{\ell\in Y_D}\Lambda_B(u+d\ell)e_D(\ell\rho)= \mathcal{F}_D\left(\Lambda_B(u+d \cdot)\right)(-\rho);
		\end{equation}
		contrast this with the function introduced in \cite[\S 6]{FKMS26}, where the authors considered
		\[ K_\theta(u) := e(u \theta/N) \cdot \mathbb{E}_{l \in Y_D} w_B(ld + u)e(l\theta/D) 
		= e(u \theta/N) \tilde{K}_{d,\theta}(u).\] 
		Consequently, the following lemma is a direct consequence of \cite[Lemma 6.3]{FKMS26}.
		\begin{lemma}\label{lem:fibre-envelope-new}
			There is \(H_d:Y_d\to[0,\infty)\) such that
			\begin{equation}\label{eq:fibre-envelope-new}
				\sup_{\rho\in Y_D}|\widetilde K_{d,\rho}(u)|\le H_d(u),
				\qquad
				\mathbb{E}_{u\in Y_d}H_d(u)\lesssim Q^{o(1)-1}.
			\end{equation}
		\end{lemma}
		
		Consider families of ``structured atoms'' 
		\begin{align}
			\mathcal{A}_d := \{  v(a) e_D(\vartheta(a) s), \; \; \; a \in \mathbb{Z}/d, \ \vartheta(a), s \in \mathbb{Z}/D, \ \left \| v \right \|_{L^2(Y_d)}= 1 \};
		\end{align}
	as in \cite[\S 5]{FKMS26}.

		The majorant $H_d$ controls the structured atoms, as the following lemma makes precise.
		
		\begin{lemma}\label{l:structured-bound}
			Suppose $\beta_d > 0$, $\Phi\in\mathcal A_d$, and $G : Y_N \to \mathbb{C}$. Then
			\begin{align}
				\| A_B(\Phi,G) \|_{L^1(Y_N)} \lesssim Q^{o(1)-1} \|G \|_{L^2(Y_N)}.
			\end{align}
		\end{lemma}
		For technical reasons, we will extend $\Phi$ and $G$ periodicially $\mod N$, and prove an estimate on the right hand side with $\| G \|_{L^2((-N,2N]} \approx \| G \|_{L^2(Y_N)}$; this will avoid overlap/periodicity issues.
		\begin{proof}
			Because $\Phi \in \mathcal{A}_d$ we have that $\Phi(a+ds)=v(a)e_D(\theta(a)s)$ with $\left \| v \right \|_{L^2(Y_d)}=1.$ Then, expanding the definition of $A_B(\Phi,G)(a'+ds')$ we get that it is controlled in absolute value by \[ \begin{split}
				& \EE_{a \in Y_d}\left |v(a'-a)e_{D}(\theta(a'-a)s')\right | \left|  \EE_{s \in Y_D}\Lambda_B(a+ds)e_D(-\theta(a'-a)s)G(a'+a+d(s'+s))\right|\\ & =\EE_{a \in Y_d}\left|v(a'-a)   \mathcal{F}_D\left(\Lambda_B(a+d \cdot)G(a'+a+d(s'+\cdot))\right)(\theta(a'-a))\right|.
			\end{split}\] 
			Taking the $L^1(Y_D)$ norm of the expression above we have that \[ \begin{split}
				&	\left \| \EE_{a \in Y_d}\left|v(a'-a)  \mathcal{F}_D\left(\Lambda_B(a+d \cdot)G(a'+a+d(s'+\cdot))\right)(\theta(a'-a))\right|\right \|_{L^1_s(Y_D)} \\ &\leq \EE_{a \in Y_d}|v(a'-a)| \EE_{s' \in Y_D}|\mathcal{F}_D\left(\Lambda_B(a+d \cdot)G(a'+a+d(s'+\cdot))\right)(\theta(a'-a))|  \\ & \leq \EE_{a \in Y_d}|v(a'-a)| \left \| \sum_{\xi \in Y_D} \mathcal{F}_D(\Lambda_B(a+d \cdot))(\theta(a'-a)-\xi) \mathcal{F}_D(G(a'+a+d \cdot))(\xi)e_{D}(\xi s') \right \|_{L^2_s(Y_D)} \\ & \leq \EE_{a \in Y_d}|v(a'-a)| \left(\sum_{\xi \in Y_D}\left|\mathcal{F}_D(\Lambda_B(a+d \cdot))(\theta(a'-a)-\xi) \mathcal{F}_D(G(a'+a+d \cdot))(\xi)\right|^2\right)^{\frac{1}{2}}  \\ & \leq \EE_{a \in Y_d}|v(a'-a)| H_d(a) \left(\EE_{s \in Y_D}|G(a'+a+ds)|^2\right)^{\frac{1}{2}}. \end{split}\] Finally, the $L^2(Y_d)$ norm of the majorant we have computed above is controlled by \[ \left \| \EE_{a \in Y_d}|v(a'-a)| H_d(a) \left(\EE_{s \in Y_D}|G(a'+a+ds)|^2\right)^{\frac{1}{2}} \right \|_{L^1_{a'}(Y_d)} \leq \left \| G \right \|_{L^2(Y_N)} \left \| v \right \|_{L^2(Y_d)} \left \| H_d \right \|_{L^1(Y_d)}   \] where the last step is due to Young's convolution inequality. The proof is completed by invoking the conclusion of Lemma \ref{lem:fibre-envelope-new}.
			
		\end{proof}

		The inverse theorem \cite[Lemma 5.5]{FKMS26} yields the following greedy decomposition, which is proved as in \cite[Proposition 5.7]{FKMS26}.
		
		\begin{proposition}[Greedy Decomposition]\label{p:greedy}
			If $F\in L^2(Y_N)$ and $t>0$, there are structured atoms, $\{ \Phi_j : j \leq J\}\subset\bigcup_{d:\,\beta_d>0}\mathcal A_d$, coefficients $\{ c_j \}$ so that
			\begin{align}
				\sum_{j \leq J} |c_j|^2 \leq \|F \|_{L^2(Y_N)}^2, \qquad J \leq t^{-2}
			\end{align}
			and a decomposition
			\begin{align}
				F = \sum_{j \leq J} c_j \Phi_j + R,\quad \text{where } \| M_d R \|_{L^2(Y_N)} \leq t \|F \|_{L^2(Y_N)}.
			\end{align}
			
		\end{proposition}
		With these preparations in mind, we are prepared to estimate $d_B(1)$:
		\begin{proposition} \label{l:d4}
			We have that 
			\begin{align}
				\| A_B(F,G) \|_{L^1(Y_N)} &\lesssim |B|^{1/2} Q^{o(1) - 1} \|F \|_{L^2(Y_N)} \| G \|_{L^2(Y_N)}, 
			\end{align}
			and in particular
			\begin{align}
				d_B(1) \leq |B|^{1/2} Q^{o(1) -1}.
			\end{align}
		\end{proposition}
		\begin{proof}
			Apply Proposition \ref{p:greedy} with $t = |B|^{-1/2}$ so that the $\ell^1$ sum over structured coefficients contributes 
			\[ \sum_{j=1}^J|c_j| \leq J^{1/2}(\sum_{j=1}^J|c_j|^2)^{1/2}\leq |B|^{1/2}\|F\|_{L^2(Y_N)};\]
			hence the principal part of $F$ comprised of the structured atoms contributes
			\begin{align}
				\left \| A_B\left (\sum_{j=1}^J c_j \Phi_j,G \right )  \right \|_{L^1(Y_N)} \leq   |B|^{1/2} Q^{o(1) -1 } \|F \|_{L^2(Y_N)} \|G \|_{L^2(Y_N)}.
			\end{align}
			Since the remainder satisfies all the fibre bounds, \eqref{eq:divisorsumdecomp} and Lemma~\ref{l:unweighted-fibre} imply
			\begin{align}
				\|A_B(R,G)\|_{L^1(Y_N)}
				&\leq\sum_{d\mid N}\beta_d\|A_d(R,G)\|_{L^1(Y_N)}\lesssim\|G\|_{L^2(Y_N)}
				\sum_{d\mid N}\beta_d\|M_dR\|_{L^2(Y_N)}\\
				&\lesssim Q^{o(1)-1}|B|^{1/2}
				\|F\|_{L^2(Y_N)}\|G\|_{L^2(Y_N)},
			\end{align}
			completing the proof.
		\end{proof}
		
		The next section will be concerned with appropriately propagating the gain furnished by our periodic analysis.
		
		\section{Transference}\label{s:transference}
		We recall the vector-valued framework of \cite{AmentaUraltsev2020,DPLMV}, which will transfer the periodic gain to the long operators in \eqref{e:longgoal}.
		\subsection{A Word From Our Sponsors: Time Frequency Analysis}\label{s:TFA} 
		
		Let $P\in\mathbb{N}$ be fixed, let $1<q_1,q_2,q_3<\infty$, and set
		\[
		Y_P:=\mathbb{Z}/P,\qquad X_i:=L^{q_i}(Y_P),
		\]
		where $Y_P$ carries normalized counting measure.
		Given a bounded form
		\[
		\Pi:X_1\times X_2\times X_3\to\mathbb{C}
		\]
		that is complex-linear in its first two arguments and conjugate-linear
		in its third, define
		\[
		\|\Pi\|:=\sup\bigl\{|\Pi(f_1,f_2,f_3)|:
		\|f_i\|_{X_i}\leq1\bigr\}.
		\]
		Because $\Pi$ is conjugate-linear in its third argument, there is a
		unique bounded bilinear map
		\[
		A_\Pi:X_1\times X_2\to L^{q_3'}(Y_P)
		\]
		such that
		\[
		\Pi(f_1,f_2,f_3)
		=\langle A_\Pi(f_1,f_2),f_3\rangle_{Y_P}
		:=\mathbb{E}_{y\in Y_P}A_\Pi(f_1,f_2)(y)\overline{f_3(y)},
		\]
		where $1/q_3+1/q_3'=1$ and $\|A_\Pi\|=\|\Pi\|$.
		
		Given Schwartz functions which map into $X_i$, i.e.\ those of the form
		\begin{align}\label{e:tensor}
			F(x;n) = \sum_{A \in \mathbb{Z}/P} f_{A}(x) e(A/P n), 
		\end{align} where $\{ f_A \}_{A \in \Z/P} $ is a collection of Schwartz functions
		and a multiplier
$	m \in C^{M}(\mathbb{R}\setminus \left\{0\right\}).$
		Define
		\begin{align}
			\Lambda_m^{\Pi}(F_1,F_2,F_3)
			:=\int_{\mathbb{R}^2} &m(\xi_2-\xi_1)\cdot \Big( \Pi\bigl(\widehat{F_1}(\xi_1;\cdot),
			\widehat{F_2}(\xi_2;\cdot),
			\widehat{F_3}(\xi_1+\xi_2;\cdot)\bigr) \Big)
			\,d\xi_1\,d\xi_2,
		\end{align}
		where the Fourier transform of \eqref{e:tensor} is defined
		\begin{align}
			\widehat{F}(\xi;n) := \sum_{A \in \mathbb{Z}/P} \widehat{f_{A}}(\xi) e(A/P n).
		\end{align}
		
		The following specialization of \cite[Theorem 1.4]{DPLMV} contains
		the Lebesgue exponents needed below, while its constant is independent
		of the size of the cyclic group.
		\begin{theorem}\label{t:finiteBV}
			If $1<s<\infty$ and $\Pi$ is bounded on
			$L^2(Y_P)\times L^2(Y_P)\times L^{s}(Y_P)$, then
			\begin{equation}\label{e:vector-221}
				\begin{split}
					|\Lambda_m^\Pi(F_1,F_2,F_3)|
					&\lesssim_s  \mathbf{S}_M(m) \cdot \|\Pi\| \cdot
					\|F_1\|_{L_x^2(\mathbb R;L_y^2(Y_P))}
					\|F_2\|_{L_x^2(\mathbb R;L_y^2(Y_P))}\\
					&\qquad\qquad\times
					\|F_3\|_{L_x^\infty(\mathbb R;L_y^{s}(Y_P))},
				\end{split}
			\end{equation}
			where $M$ is sufficiently large and
			$\mathbf{S}_M(m):=\max_{0\leq j\leq M}
			\sup_{t\neq0}|t|^j|\partial^jm(t)|$.
		\end{theorem}
		\begin{proof}
			The complex-trilinear form
			\[
			(u_1,u_2,u_3)\longmapsto\Pi(u_1,u_2,\overline{u_3})
			\]
			has norm $\|\Pi\|$, and the identity
			\[
			\widehat{\overline{F_3}}(-\xi_1-\xi_2;\cdot)
			=\overline{\widehat{F_3}(\xi_1+\xi_2;\cdot)}
			\]
			converts the complex-trilinear convention of \cite{DPLMV} into ours
			when the third test function is $\overline{F_3}$.
			The first two spaces are Hilbert spaces, while $L^{s}(Y_P)$ is
			$R$-intermediate UMD for every $R>\max\{s,s'\}$, with uniform constants
			in $P$.
			Taking the intermediate exponents $(2,2,R)$ in
			\cite[Theorem 1.4]{DPLMV}, the associated parameter is $\rho=1/R$,
			and the four vertices
			\[
			(-1/2,1/2),\quad(1/2,-1/2),\quad
			(1/2+1/R,1/2),\quad(1/2,1/2+1/R)
			\]
			already enclose $(1/2,1/2)$ in their convex hull's interior.
			The conclusion of \cite[Theorem 1.4]{DPLMV} for the corresponding bilinear operator applies at Lebesgue exponents
			$(2,2,1)$, which yields \eqref{e:vector-221} after duality against
			$L_x^\infty L_y^{s'}$ and normalization of $m$ and $\Pi$.
		\end{proof}

		We now address aliasing issues;
		
		\medskip 
		
		To do so, choose finitely many smooth functions $\chi_{\ell}$ on $\mathbb R$, each
		supported on an interval $I_\ell$ of length less than $10^{-3}$, such that
		\begin{align}
			\sum_{\ell\leq M_0}\sum_{j\in\mathbb Z}\chi_{\ell}(\theta-j)=1,
			\qquad \theta\in\mathbb R.
		\end{align}
		In particular, the periodizations of the $\chi_{\ell}$ form a partition of unity on
		$\mathbb T$. For a positive integer $P$, define
		\begin{align}
			\widehat{P_{\ell,P}h}(\beta)
			:=\Big(\sum_{j\in\mathbb Z}\chi_{\ell}(P\beta-j)\Big)
			\widehat h(\beta),\qquad \beta\in\mathbb T,
		\end{align}
		and for a finitely supported sequence $h$, let $h_A$ be the Schwartz function
		on $\mathbb R$ with
		\begin{align}\label{e:chart-def}
			\widehat{h_A}(\alpha)
			:=\chi_{\ell}(P\alpha)\widehat h(\alpha+A/P),
			\qquad A\in Y_P,
		\end{align}
		where $\widehat h$ on the right is extended periodically to $\mathbb R$;
		Fourier inversion yields
		\begin{align}
			P_{\ell,P}h(n)=\sum_{A\in Y_P}h_A(n)e(An/P),
			\qquad n\in\mathbb Z.
		\end{align}
		The functions we will be most interested in are the lifts
		\begin{align}
			H(x;y):=\sum_{A\in Y_P}h_A(x)e(Ay/P),
			\qquad x\in\mathbb R,\quad y\in Y_P,
		\end{align}
		whose norms we will control by $h$ below; we normalize all averages in $y$, so that the following bounds are
		uniform in the cyclic group in both their discrete and continuous forms.
		
		\begin{lemma}\label{l:lift}
			Regarding $\ell$ as fixed, for each $ p \in [1,\infty]$ we have
			\begin{align}
				\|H\|_{\ell^p_x(\mathbb Z;L^p_y(Y_P))}
				+\|H\|_{L^p_x(\mathbb R;L^p_y(Y_P))}
				\lesssim_p\|h\|_{\ell^p(\mathbb Z)},
			\end{align}
			with a constant independent of $P$.
		\end{lemma}
		\begin{proof}
			Writing $\kappa_\ell:=\chi_{\ell}^\vee$, Fourier inversion yields
			\[
			H(x;y)=\sum_{\substack{n\in\mathbb Z\\n\equiv y\mod P}}
			h(n)\kappa_\ell\big((x-n)/P\big),
			\]
			where the factor $P^{-1}$ from dilation cancels the character sum.
			Since $\kappa_\ell$ is Schwartz,
			\[
			C_\ell:=\sup_{t\in\mathbb R}\sum_{j\in\mathbb Z}
			|\kappa_\ell(t-j)|<\infty,
			\]
			and weighted H\"older implies, for $1\leq p<\infty$,
			\[
			|H(x;y)|^p\leq C_\ell^{p-1}
			\sum_{n\equiv y\mod P}|h(n)|^p
			|\kappa_\ell((x-n)/P)|.
			\]
			Averaging in $y$ and integrating in $x$, we obtain
			\[
			\|H\|_{L_x^pL_y^p}^p
			\leq C_\ell^{p-1}\|\kappa_\ell\|_{L^1}\|h\|_{\ell^p}^p,
			\]
			because the average contributes $P^{-1}$ and the integral contributes $P$.
			The discrete estimate follows by the same argument, since
			$P^{-1}\sum_{m\in\mathbb Z}|\kappa_\ell(m/P)|$ is uniformly bounded.
			For $p=\infty$, the bound
			$|H(x;y)|\leq C_\ell\|h\|_{\ell^\infty}$ follows from the first
			formula, which completes both assertions.
		\end{proof}

		Since $Y_P$ is a probability space, the preceding lemma also controls
		the smaller inner exponents that occur in \eqref{e:vector-221}.
		\begin{corollary}\label{c:lift}
			If $1\leq s\leq p$ and $p \in [1,\infty]$, then
			\begin{align}
				\|H\|_{L^p_x(\mathbb R;L^s_y(Y_P))}
				\lesssim_p\|h\|_{\ell^p(\mathbb Z)}.
			\end{align}
		\end{corollary}
		
		\subsection{Exact transference of a periodic block}
		We retain $\tau$ and $\epsilon$ from Section~\ref{s:red}, while for a
		fixed block $B$ we abbreviate \[ N=\mathcal Q_B, \qquad K_Q=K_{Q;\tau}, \qquad \mathcal H_{Q,B}=\mathcal H_{Q;\tau,B}.\]
		Because $\psi$ is smooth, compactly supported away from zero, and odd,
		the dyadic multiplier estimates and the product rule imply
		\begin{align}\label{e:kernel-BV}
			\mathbf{S}_M(\widehat K_Q)&\lesssim_M1, \qquad \operatorname{supp}\widehat K_Q
			\subset[-c/N(Q),c/N(Q)],
		\end{align}
		uniformly in $Q$ and $J$, where $c$ depends only on
		$\operatorname{supp}\widehat\varphi$.
		As $\kappa=\tau-\epsilon$, we have
		\[
		\log N\lesssim Q^{\tau-\epsilon}\log Q,
		\qquad \log N(Q)=Q^{\tau-\epsilon/10},
		\]
		so that $N/N(Q)\to0$ rapidly, and we may assume $cN/N(Q)<10^{-3}$ after
		absorbing finitely many values of $Q$ using \eqref{e:progression}.
		
		Fix a triple $(\ell_1,\ell_2,\ell_3)$, and let
		$f_a,g_b,h_c$ denote the scalar coefficients from \eqref{e:chart-def},
		with $a,b,c\in Y_N$ represented by integers in $\{0,\ldots,N-1\}$.
		We write
		\begin{align}\label{e:QBform}
			\mathcal T_{\ell_1,\ell_2,\ell_3}
			:=\big\langle\mathcal H_{Q,B}(P_{\ell_1,N}f,P_{\ell_2,N}g),
			P_{\ell_3,N}h\big\rangle_{\ell^2},
		\end{align}
		where the pairing is the unnormalized sum against the complex conjugate.
		If
		\[
		M_B(\theta):=\sum_nK_Q(n)\Lambda_B(n)e(-n\theta),
		\]
		then Fourier inversion expresses this form as
		\begin{equation}\label{e:chart-expanded}
			\begin{split}
				\mathcal T_{\ell_1,\ell_2,\ell_3}
				=\sum_{a,b,c\in Y_N}\sum_{k\in\mathbb Z}
				\int_{\mathbb R^2}&\widehat f_a(\alpha)\widehat g_b(\beta)
				\overline{\widehat h_c\big(\alpha+\beta+(a+b-c)/N-k\big)}\\
				&\qquad\times M_B\big((a-b)/N+\alpha-\beta\big)
				\,d\alpha\,d\beta,
			\end{split}
		\end{equation}
		in which the sum over $k$ records the periodization of the output frequency.
		Since $K_Q$ is Schwartz for finite $J$, Poisson summation is applicable
		and yields
		\begin{align}\label{e:sampled-multiplier}
			M_B(\theta)
			=\sum_{\ell\in Y_N}\mathcal F_N(\Lambda_B)(\ell)
			\sum_{j\in\mathbb Z}\widehat K_Q(\theta-\ell/N-j);
		\end{align}
		the support condition in \eqref{e:kernel-BV} restricts a surviving term to
		\[
		|N(\alpha-\beta)-d|\leq cN/N(Q),\qquad
		d:=\ell-a+b+jN\in\mathbb Z,
		\]
		so there is at most one possible $d$ for the entire pair $(\ell_1,\ell_2)$, because $I_{\ell_1}-I_{\ell_2}$ has length less than $2\cdot10^{-3}$. If such an integer exists, the multiplier reduces exactly to
		\begin{align}\label{e:one-alias}
			M_B\big((a-b)/N+\alpha-\beta\big)
			=\mathcal F_N(\Lambda_B)(a-b+d)
			\widehat K_Q(\alpha-\beta-d/N),
		\end{align}
		on the supports of $\widehat f_a$ and $\widehat g_b$, where the argument of $\mathcal F_N$ is taken modulo $N$. In particular, whenever the integrand in \eqref{e:chart-expanded} is nonzero, the chart definitions imply
		\[
		N\alpha\in I_{\ell_1},\qquad
		N\beta\in I_{\ell_2},\qquad
		N\zeta\in I_{\ell_3},
		\qquad
		\zeta=\alpha+\beta+\frac{a+b-c}{N}-k,
		\]
		so
		\begin{align}\label{e:int00}
			v:=N(\alpha+\beta-\zeta)=Nk-a-b+c
			\in \mathbb{Z}\cap
			\bigl(I_{\ell_1}+I_{\ell_2}-I_{\ell_3}\bigr),
		\end{align}
		where the interval on the right has length less than
		$3\cdot 10^{-3}<1$, so in particular the set on the right of \eqref{e:int00} contains at most one point; since this interval depends only on the chart triple, if this set is non-empty, all nonzero contributions have the same value of $v$ (otherwise the chart contribution vanishes). In particular, in the non-trivial case the surviving indices satisfy $a+b-c+v\equiv0\mod N$,
		and \eqref{e:chart-expanded} becomes
		\begin{equation}\label{e:exact-chart-integral}
			\begin{split}
				\mathcal T_{\ell_1,\ell_2,\ell_3}
				=\sum_{\substack{a,b,c\in Y_N\\a+b-c+v\equiv0\mod N}}
				&\mathcal F_N(\Lambda_B)(a-b+d)\int_{\mathbb R^2}
				\widehat K_Q(\alpha-\beta-d/N)\\
				&\quad\times\widehat f_a(\alpha)\widehat g_b(\beta)
				\overline{\widehat h_c(\alpha+\beta-v/N)}
				\,d\alpha\,d\beta.
			\end{split}
		\end{equation}
		
		To identify this expression with a vector-valued form, define
		\begin{equation}\label{e:twisted-Pi}
			\Pi_{B,d,v}(F_1,F_2,F_3)
			:=\mathbb E_{y,t\in Y_N}\Lambda_B(t)e_N(vy-dt)
			F_1(y-t)F_2(y+t)\overline{F_3(y)}.
		\end{equation}
		Writing $e_a(y):=e_N(ay)$, we have
		\[
		\Pi_{B,d,v}(e_a,e_b,e_c)
		=\mathbf1_{a+b-c+v\equiv0\mod N} \cdot
		\mathcal F_N(\Lambda_B)(a-b+d),
		\]
		and the identity
		\[
		\begin{aligned}
			\Pi_{B,d,v}(F_1,F_2,F_3)
			&=\mathbb E_y A_B(e_dF_1,F_2)(y)e_{v-d}(y)\overline{F_3(y)}\\
			&=\Pi_B(e_dF_1,F_2,e_{d-v}F_3)
		\end{aligned}
		\]
		shows that, for $1<s\leq2$,
		\begin{align}\label{e:twisted-Pi-norm}
			\|\Pi_{B,d,v}\|_{L^2(Y_N)\times L^2(Y_N)\times L^{s'}(Y_N)}
			\leq d_B(s).
		\end{align}
		Character multiplication preserves the relevant norms, so this estimate is uniform in the two aliasing integers, $d,v$; the two integers have distinct roles, since $d$ comes from sampling the kernel and $v$ comes from periodizing the output frequency.
		
		Let $F,G,H$ be the lifts of the three chart inputs, and set
		\[
		\begin{aligned}
			m(u)&:=\widehat K_Q(-u),\\
			F_{d,N}(x;y)&:=e(-dx/N)F(x;y),\qquad
			H_{d,v,N}(x;y):=e((v-d)x/N)H(x;y).
		\end{aligned}
		\]
		Making the change of variables $\alpha=\alpha'+d/N$ in
		\eqref{e:exact-chart-integral}, we obtain the exact identity
		\begin{align}\label{e:chart-identification}
			\mathcal T_{\ell_1,\ell_2,\ell_3}
			=\Lambda_m^{\Pi_{B,d,v}}(F_{d,N},G,H_{d,v,N}),
		\end{align}
		since
		\[
		\begin{aligned}
			\widehat{F_{d,N}}(\alpha;y)
			&=\sum_a\widehat f_a(\alpha+d/N)e_a(y), \qquad \text{ and} \\
			\widehat{H_{d,v,N}}(\alpha+\beta;y)
			&=\sum_c\widehat h_c(\alpha+\beta+(d-v)/N)e_c(y).
		\end{aligned}
		\]
		Because $\Pi_{B,d,v}$ is conjugate-linear in its third argument, its
		expansion conjugates the coefficients in the second identity, producing
		the factor $\overline{\widehat h_c(\alpha+\beta+(d-v)/N)}$ required by
		\eqref{e:exact-chart-integral} after the change of variables.

		Applying Theorem~\ref{t:finiteBV} and Corollary~\ref{c:lift}, while
		using the fact that the displayed modulations preserve mixed norms,
		we conclude that
		\begin{align*}
			|\mathcal T_{\ell_1,\ell_2,\ell_3}|
			&\lesssim_s d_B(s)
			\|F\|_{L_x^2L_y^2}\|G\|_{L_x^2L_y^2}
			\|H\|_{L_x^\infty L_y^{s'}}\\
			&\lesssim_s d_B(s)\|f\|_{\ell^2}\|g\|_{\ell^2}\|h\|_{\ell^\infty}.
		\end{align*}
		Since the number of chart triples is fixed, summation and duality imply
		\begin{align}\label{e:block-transfer}
			\|\mathcal H_{Q;\tau,B}(f,g)\|_{\ell^1}
			\lesssim_s d_B(s)\|f\|_{\ell^2}\|g\|_{\ell^2} \lesssim_s |B|^{1/2} Q^{o(1)-1/s} \|f\|_{\ell^2}\|g\|_{\ell^2} 
		\end{align}
		with the final inequality deriving from Proposition~\ref{p:interpsetup}.
		
		\begin{proof}[Proof of Proposition~\ref{p:key0}]
			Since the partition has $O(Q^{1-\kappa})$ blocks of cardinality
			$O(Q^\kappa)$, \eqref{e:block-transfer} implies
			\begin{align}\label{e:Q-saving}
				\|\mathcal H_{Q;\tau}(f,g)\|_{\ell^1}
				&\lesssim_s Q^{1-\kappa/2-1/s+o(1)}
				\|f\|_{\ell^2}\|g\|_{\ell^2}\notag\\
				&\lesssim_s Q^{-\tau/2+\epsilon/2+(1-1/s)+o(1)}
				\|f\|_{\ell^2}\|g\|_{\ell^2} \\
				& \lesssim_s Q^{-\tau/2+o(1)} \|f\|_{\ell^2}\|g\|_{\ell^2}
			\end{align}
			for $s$ sufficiently close to $1$ and $\epsilon >0$ sufficiently small.
			The unsmoothed estimate follows from \eqref{e:smootherror}, since
			$B^*$ maps $\ell^2\times\ell^2$ to $\ell^1$.
		\end{proof}
		We are now ready to prove Theorem \ref{t:maintheorem}.
		\section{Proof of Theorem~\ref{t:maintheorem}}\label{s:final}
		\begin{proof}[Proof of Theorem~\ref{t:maintheorem}]
			Fix a triple $(p,q,r)\in\mathcal R$, and retain a finite upper truncation $J$. We will obtain independent of $J$  bounds for   $ \| H^J \|_{\ell^p \times \ell^q \to \ell^r}.$ 
			The error and bounded-complexity terms in \eqref{e:exactred} are
			controlled by \eqref{e:error-summable} and \eqref{e:progression}, so it
			remains to obtain summable bounds for the two high-complexity components.
			
			\emph{The case $r\geq1$.}
			Taking $\tau_0=1/6$ and summing at most $Q^{1/6}$ single scales,
			\eqref{e:single-L1} and \eqref{e:long-unsmoothed}  control the contributions of $H^J_{Q,\mathrm{moderate};1/6}$ and $U_{Q;1/6}$ respectively:
			\[
			\|H^J_{Q,\mathrm{moderate};1/6}(f,g)\|_{\ell^1}
			+\|U_{Q;1/6}(f,g)\|_{\ell^1}
			\lesssim Q^{-1/20}\|f\|_{\ell^2}\|g\|_{\ell^2}.
			\]
			Since $1<p,q<\infty$, every target reciprocal tuple is a convex
			combination of $(1/2,1/2,1)$ and another Banach tuple, so interpolation with \eqref{e:progression} yields the desired power savings in $Q$.
			
			\emph{The case $r<1$.}
			By symmetry we assume $p\leq q$, and put $t:=1/r-1$.
			The defining inequalities of $\mathcal R$ imply
			\begin{align}\label{e:range-gaps}
				1-\frac1p>2t,\qquad 0<t<\frac15,
			\end{align}
			so we may choose $\theta$ satisfying
			\[
			\frac83t<\theta<\min\left\{\frac{2}{p'},\,1-2t\right\}.
			\]
			Define $(p_0,q_0,r_0)$ by
			\[
			\frac1{p_0}=\frac{1/p-\theta/2}{1-\theta},\qquad
			\frac1{q_0}=\frac{1/q-\theta/2}{1-\theta},\qquad
			\frac1{r_0}=\frac{1/r-\theta}{1-\theta},
			\]
			so that $p_0,q_0>1$, $2/3<r_0<1$, and interpolation with
			$(2,2,1)$ at weight $\theta$ reaches $(p,q,r)$.
			With $\tau=3/4$, interpolation of \eqref{e:long-unsmoothed} and
			\eqref{e:progression} implies
			\begin{align}\label{eq:singlecomplexitylongquasibanach}
				\|U_{Q;\tau}(f,g)\|_{\ell^r}
				\lesssim Q^{t-3\theta/8+o(1)}\|f\|_{\ell^p}\|g\|_{\ell^q},
			\end{align}
			whose exponent is negative by the choice of $\theta$.
			The same interpolation between \eqref{e:single-L1} and
			\eqref{e:single-rough} yields a single-scale norm
			$O(Q^{-\theta/4+o(1)})$, so for any $0<\tau_0<r\theta/4$ we have
			\begin{align}\label{e:negpower}
				\|H^J_{Q,\mathrm{moderate};\tau_0}(f,g)\|_{\ell^r}
				\lesssim Q^{\tau_0/r-\theta/4+o(1)}
				\|f\|_{\ell^p}\|g\|_{\ell^q}.
			\end{align}
			After increasing $\mathbf C$ we may assume $b_Q<Q^{\tau_0}$,
			while this initial group already has summable decay.
			
			We now delicately aggregate scales for interpolation purposes; this move is not needed in the case $r \geq 1$. Specifically, choose a finite increasing sequence
			$\tau_0<\tau_1<\cdots<\tau_M=3/4$ such that
			\begin{align}\label{e:group-cutoffs}
				2\tau_jt<\frac{\tau_{j-1}}{p'},
				\qquad 1\leq j\leq M,
			\end{align}
			which is possible by \eqref{e:range-gaps}, with $M$ depending only
			on the target tuple.
			For each $j$ put
			\[
			V_{Q,j}^J(f,g):=
			\sum_{\substack{Q^{\tau_{j-1}}\leq k<Q^{\tau_j}\\k\leq J}}
			\sum_n\psi_k(n)w_Q(n)f(x-n)g(x+n),
			\]
			so that the moderate contribution is the sum of the initial group
			and these $M$ consecutive groups.
			For $0<\delta<\min\{1/2,1-1/p\}$, the three estimates we need are
			\begin{align}
				\|V_{Q,j}^J(f,g)\|_{\ell^1}
				&\lesssim Q^{-\tau_{j-1}/2+o(1)}
				\|f\|_{\ell^2}\|g\|_{\ell^2},\label{e:secondinterp}\\
				\|V_{Q,j}^J(f,g)\|_{\ell^1}
				&\lesssim_\delta Q^{o(1)}
				\|f\|_{\ell^{1/(1-\delta)}}\|g\|_{\ell^{1/\delta}},
				\label{e:thirdinterp}\\
				\|V_{Q,j}^J(f,g)\|_{\ell^{1/(2-2\delta)}}
				&\lesssim_\delta Q^{2\tau_j(1-\delta)+o(1)}
				\|f\|_{\ell^{1/(1-\delta)}}
				\|g\|_{\ell^{1/(1-\delta)}}.\label{e:fourthinterp}
			\end{align}
			The first follows from \eqref{e:long-unsmoothed} with upper truncation
			$\min\{J,\lceil Q^{\tau_j}\rceil-1\}$, while the second is
			\eqref{e:progression} and the third follows by the quasi-triangle
			inequality from \eqref{e:single-rough}.
			
			Interpolate these estimates with weights
			\[
			\lambda_0:=\frac{2(1/p'-\delta)}{1-2\delta},\qquad
			\lambda_1:=\frac{1/p-1/q}{1-2\delta},\qquad
			\lambda_2:=\frac{t}{1-2\delta},
			\]
			which are nonnegative, sum to one, and satisfy
			\[
			\left(\frac1p,\frac1q,\frac1r\right)
			=\lambda_0\left(\frac12,\frac12,1\right)
			+\lambda_1(1-\delta,\delta,1)
			+\lambda_2(1-\delta,1-\delta,2-2\delta).
			\]
			The resulting bound is
			\begin{align}\label{e:group-decay}
				\|V_{Q,j}^J(f,g)\|_{\ell^r}
				&\lesssim_{p,q,\delta}Q^{\sigma_j(\delta)+o(1)}
				\|f\|_{\ell^p}\|g\|_{\ell^q},\notag\\
				\sigma_j(\delta)
				&:=\frac{-\tau_{j-1}(1/p'-\delta)
					+2\tau_j(1-\delta)t}{1-2\delta}.
			\end{align}
			Since \eqref{e:group-cutoffs} ensures $\lim_{\delta \downarrow0} \sigma_j(\delta)<0$ for every $j \leq M$,
			we choose one sufficiently small $\delta>0$ for which all these
			exponents remain negative.
			The number of groups is fixed, so the identity
			\[
			H^J_{Q,\mathrm{moderate};3/4}
			=H^J_{Q,\mathrm{moderate};\tau_0}+\sum_{j=1}^M V_{Q,j}^J
			\]
			and the $r$-quasi-triangle inequality combine
			\eqref{eq:singlecomplexitylongquasibanach}, \eqref{e:negpower}, and
			\eqref{e:group-decay} into
			\begin{align}\label{e:complexity-decay}
				\|H^J_{Q,\mathrm{moderate};3/4}(f,g)\|_{\ell^r}
				+\|U_{Q;3/4}(f,g)\|_{\ell^r}
				\lesssim_{p,q}Q^{-\sigma}\|f\|_{\ell^p}\|g\|_{\ell^q}
			\end{align}
			for some $\sigma>0$, after choosing every power loss smaller than
			the finitely many strict gaps.
			
			In either range, summation in \eqref{e:exactred}, with the
			$\min\{r,1\}$-triangle inequality, now proves
			\[
			\sup_{J\geq1}\|H^J(f,g)\|_{\ell^r}
			\lesssim_{p,q,r}\|f\|_{\ell^p}\|g\|_{\ell^q},
			\]
			since the errors decay like $k^{-100}$ and the high-complexity
			terms have a strictly negative power of $Q$.
			For finitely supported inputs, the identity
			$\sum_{k\geq1}\psi_k(n)=1/n$ identifies the limit with
			$\mathcal H_{\mathbb P}(f,g)$, and density extends the resulting
			bounded bilinear operator to $\ell^p\times\ell^q$.
		\end{proof}
		
		\subsection{AI Declaration}
		Chat GPT 6.0 Astra assisted with the above interpolation argument, used to extend the range of exponents to the full convex hull of $\mathcal{R}$.
		
		\bibliography{primeBHTbib}
		\bibliographystyle{amsplain}

	\end{document}